\documentclass{amsart}
\pdfoutput=1
\usepackage{amssymb}
\usepackage{thmtools, thm-restate}
\usepackage{graphicx}
\usepackage{hyperref} 

\declaretheorem[numberwithin=section]{theorem}

\declaretheorem[numberlike=theorem]{conjecture}
\declaretheorem[numberlike=theorem, style=definition]{definition}

\declaretheorem[numberlike=theorem, style=remark]{remark}

\numberwithin{equation}{section}

\title{Monotone Hurwitz Numbers and the HCIZ Integral}
\date{\today}

\author{I. P. Goulden}
\address{Department of Combinatorics \& Optimization. University of Waterloo, Canada}
\email{ipgoulden@uwaterloo.ca}

\author{M. Guay-Paquet}
\address{Department of Combinatorics \& Optimization. University of Waterloo, Canada}
\email{mguaypaq@uwaterloo.ca}
\thanks{IPG and MG-P supported by NSERC}

\author{J. Novak}
\address{Department of Mathematics. Massachusetts Institute of Technology, USA}
\email{jnovak@math.mit.edu}

\date{\today}

\newcommand{\field}[1]{\mathbb{#1}}

\newcommand{\hur}{H}
\newcommand{\Hur}{\mathbf{H}}
\newcommand{\mon}{\vec{H}}
\newcommand{\Mon}{\vec{\mathbf{H}}}

\newcommand{\Tr}{\operatorname{Tr}}
\newcommand{\Wg}{\operatorname{Wg}}

\newcommand{\id}{\operatorname{id}}

\newcommand{\Stab}{\operatorname{Stab}}
\newcommand{\Fix}{\operatorname{Fix}}

\begin{document}

\begin{abstract}
	In this article, we prove that the complex convergence of the HCIZ free
	energy is equivalent to the non-vanishing of the HCIZ integral in a 
	neighbourhood of $z=0$.  Our approach is based on a combinatorial
	model for the Maclaurin coefficients of the HCIZ integral
	together with classical complex-analytic techniques.
\end{abstract}

\maketitle

\setcounter{tocdepth}{2}
\tableofcontents

\setcounter{section}{-1}

\section{Introduction}
The Harish-Chandra-Itzykson-Zuber integral,

	\begin{equation}
		\label{eqn:HCIZ}
		I_N(z) = \int\limits_{U(N)} e^{zN\Tr(A_NUB_NU^{-1})}\mathrm{d}U,
	\end{equation}

\noindent
is a ubiquitous special function which plays a key role in 
random matrix theory and related areas.  It enters into both
the fine-scale spectral analysis of a single random matrix and 
the macroscopic analysis of several
interacting random matrices.  For further information on these
and other aspects of the HCIZ integral, we refer the reader to the surveys
\cite{EY,Guionnet:PS,ZZ}.

In \eqref{eqn:HCIZ}, the integration is over the group of $N \times N$
unitary matrices against the normalized Haar measure, $z$ is a 
complex parameter, and $A_N, B_N$ are any two $N \times N$ complex diagonal
matrices.  Since $U(N)$ is compact, $I_N(z)$ is an entire function of the complex variable $z$.
Let $\Omega_N$ be a simply connected open set in $\field{C}$ which contains
the origin and avoids the zeros of $I_N(z)$.  The existence of such a domain is guaranteed
by the discreteness zeros of analytic functions together with the evaluation $I_N(0)=1$.  Let
$\log I_N(z)$ denote the principal branch of the logarithm of $I_N(z)$ on $\Omega_N$.
The holomorphic function 

	\begin{equation*}
		F_N(z) = \frac{1}{N^2} \log I_N(z), \quad z \in \Omega_N,
	\end{equation*}
	
\noindent
is known in the random matrix literature as the \emph{free energy} of the
HCIZ integral.

A longstanding conjecture asserts that, if $(A_N)_{N=1}^{\infty}$ and
$(B_N)_{N=1}^{\infty}$ are two sequences of complex diagonal matrices
which grow in a suitably regular fashion, then $F_N(z)$ converges 
uniformly on compact subsets of a complex neighbourhood of $z=0$.
This conjecture is implicit in the classic paper \cite{IZ}, and explicitly 
stated in \cite{Collins,CGM}.  Clearly, a necessary condition for this
conjecture to hold is the existence of a domain $\Omega$ containing
the origin such that $I_N(z)$ is non-vanishing on $\Omega$ for all but
finitely many $N$.  In this paper, we will prove that this condition is
also sufficient.

\section{The leading derivatives theorem}
Let $S(d)$ denote the symmetric group acting on $\{1,\dots,d\}$,
and identify $S(d)$ with its (right) Cayley graph as generated by
the full conjugacy class of transpositions.
Define an edge labelling of the Cayley graph as follows:
each edge corresponding to the transposition $\tau=(s\ t)$
is marked by $t$, the larger of the two numbers interchanged.
Thus, emanating from each vertex of the Cayley graph,
there is one $2$-edge, two $3$-edges, three $4$-edges, etc.

\begin{definition}
	\label{def:MonotoneWalk}
	A walk on the Cayley graph of $S(d)$ is said to be 
	\emph{monotone} if the labels of the edges it traverses
	form a weakly increasing sequence.  
\end{definition}

Given two permutations 
$\rho,\sigma \in S(d)$, we denote by $\vec{w}^r(\rho,\sigma)$
the number of $r$-step monotone walks from $\rho$ to $\sigma$.
Given two Young diagrams $\alpha,\beta \vdash d$, we denote
by 
	
	\begin{equation}
		\label{eqn:TotalCount}
		\vec{W}^r(\alpha,\beta) = \sum_{\rho \in C_\alpha} \sum_{\sigma \in C_\beta} 
		\vec{w}^r(\rho,\sigma)
	\end{equation}
		
\noindent
the number of monotone walks beginning in the conjugacy class
$C_\alpha$ and ending in the conjugacy class $C_\beta$.

	\begin{theorem}[Leading Derivatives Theorem]
	\label{thm:LeadingDerivatives}
		For any $N \geq 1$ and any $1 \leq d \leq N$, we have
		
			\begin{equation*}
				I_N^{(d)}(0) = \sum_{r=0}^{\infty} \bigg{(}  -\frac{1}{N} \bigg{)}^r
				\sum_{\alpha,\beta \vdash d} p_\alpha(A_N) p_\beta(B_N) 
				\vec{W}^r(\alpha,\beta),
			\end{equation*}
			
		\noindent
		and this series is absolutely convergent.  Equivalently, we have
		
			\begin{equation*}
				I_N(z) = 1 + \sum_{d=1}^N \frac{z^d}{d!} \sum_{r=0}^{\infty} \bigg{(}  -\frac{1}{N} \bigg{)}^r
				\sum_{\alpha,\beta \vdash d} p_\alpha(A_N) p_\beta(B_N) 
				\vec{W}^r(\alpha,\beta) + O(z^{N+1}),
			\end{equation*}
			
		\noindent
		where the $O$-term is uniform on compact subsets of $\field{C}$.
		Here

			\begin{equation*}
				p_\alpha(A_N) = \prod_{i=1}^{\ell(\alpha)} \Tr(A_N^{\alpha_i}), 
				\quad 
				p_\beta(B_N) = \prod_{j=1}^{\ell(\beta)} \Tr(B_N^{\beta_j}),
			\end{equation*}
	
		\noindent
		are the power sum symmetric functions $p_\alpha,p_\beta$ evaluated on 
		the eigenvalues of $A_N,B_N$.
	\end{theorem}

In the remainder of this section, we give the proof of Theorem \ref{thm:LeadingDerivatives}.

	\subsection{Differentiation under the integral sign}
	The derivatives of the entire function \eqref{eqn:HCIZ}
	may be computed by differentiation under the integral sign.		
	In particular, the Maclaurin coefficients of $I_N(z)$ are given by
	
		\begin{align*}
			I_N^{(d)}(0) &= N^d\int\limits_{U(N)} (\Tr A_NUB_NU^{-1})^d  \mathrm{d}U \\
			&= N^d\sum_{i,j} a_{i(1)} \dots a_{i(d)} 
			b_{j(1)} \dots b_{j(d)} \int\limits_{U(N)} |u_{i(1)j(1)} \dots u_{i(d)j(d)}|^2 \mathrm{d}U,
		\end{align*}
		
	\noindent
	where the summation is over all $N^{2d}$ pairs of functions 
	
		\begin{equation*}
			i,j:\{1,\dots,d\} \rightarrow \{1,\dots,N\},
		\end{equation*}
		
	\noindent
	and $A_N=\operatorname{diag}(a_1,\dots,a_N)$, $B_N=\operatorname{diag}(b_1,\dots,b_N)$.
				
	\subsection{The Weingarten function}
	The integration of monomial functions of matrix elements over $U(N)$ 
	can be addressed using the \emph{Weingarten convolution formula}
	of Collins and \'Sniady~\cite{CS}: 
	
		\begin{equation*}
			\int\limits_{U(N)} u_{i(1)j(1)} \dots u_{i(d)j(d)}
			\overline{u_{i'(1)j'(1)} \dots u_{i'(d)j'(d)}}\mathrm{d}U
			=\sum_{\rho,\sigma \in S(d)} \delta_{i, i'\rho} \delta_{j, j'\sigma} 
			\Wg_N(\rho^{-1}\sigma),
		\end{equation*}
		
	\noindent
	where $\Wg_N:S(d) \rightarrow \field{Q}$ is the \emph{Weingarten function}, which is given by
		
		\begin{equation*}
			\Wg_N(\pi) = \int\limits_{U(N)} u_{11} \dots u_{dd} \overline{u_{1\pi(1)} \dots u_{d\pi(d)}} \mathrm{d}U
		\end{equation*}
		
	\noindent
	for $N \geq d$.
	
	Given a function $i:\{1,\dots,d\} \rightarrow \{1,\dots,N\}$, we denote by
	$\Stab(i)$ the set of permutations $\pi \in S(d)$ such that $i\pi=i$,	
	and given a permutation $\pi \in S(d)$ we denote by 
	$\Fix(\pi)$ the set of functions
	$i:\{1,\dots,d\} \rightarrow \{1,\dots,N\}$ such that $i\pi=i$.
	Applying the Weingarten convolution formula to $I_N^{(d)}(0)$ with 
	$1 \leq d \leq N$, we obtain
	
		\begin{equation}
			\label{eqn:MaclaurinWeingarten}
			\begin{split}
			I_N^{(d)}(0) &= N^d\sum_i \sum_j a_{i(1)} \dots a_{i(d)} b_{j(1)} \dots b_{j(d)} \sum_{\rho \in \Stab(i)} \sum_{\sigma \in \Stab(j)} \Wg_N(\rho^{-1}\sigma) \\
			&=N^d \sum_\rho \sum_\sigma \Wg_N(\rho^{-1}\sigma) \sum_{i \in\Fix(\rho)} \sum_{j \in \Fix(\sigma)}  
			a_{i(1)} \dots a_{i(d)} b_{j(1)} \dots b_{j(d)}\\
			&= N^d\sum_{\rho,\sigma} \Wg_N(\rho^{-1}\sigma) p_{t(\rho)}(A_N) p_{t(\sigma)}(B_N).
			\end{split}
		\end{equation}
			
	\subsection{$1/N$-expansion of the Weingarten function}
	In \cite{Novak:BCP} (see also \cite{MN}), it was shown that, for any $N \geq d$
	and any $\pi \in S(d)$,
	the Weingarten function admits the following absolutely convergent expansion
	in powers of $1/N$:
	
		\begin{equation*}
			\Wg_N(\pi) = \frac{1}{N^d} \sum_{r=0}^{\infty} (-1)^r \frac{\vec{w}^r(\id,\pi)}{N^r}.
		\end{equation*}
		
	\noindent
	Thus
		
		\begin{align*}
			\Wg_N(\rho^{-1}\sigma)  = \frac{1}{N^d} \sum_{r=0}^{\infty} (-1)^r \frac{\vec{w}^r(\id,\rho^{-1}\sigma)}{N^r}
			= \frac{1}{N^d} \sum_{r=0}^{\infty} (-1)^r \frac{\vec{w}^r(\rho,\sigma)}{N^r}.
		\end{align*}
		
	\noindent
	Plugging this expansion into \eqref{eqn:MaclaurinWeingarten}
	and changing order of summation, we arrive at 
		
		\begin{align*}
			I_N^{(d)}(0) &= \sum_{\rho,\sigma \in S(d)} p_{t(\rho)}(A_N) p_{t(\sigma)}(B_N) \sum_{r=0}^\infty (-1)^r  \frac{\vec{w}^r(\rho,\sigma)}{N^r} \\
			&= \sum_{r=0}^\infty \bigg{(} -\frac{1}{N} \bigg{)}^r \sum_{\rho,\sigma \in S(d)}  p_{t(\rho)}(A_N) p_{t(\sigma)}(B_N) \vec{w}^r(\rho,\sigma),
		\end{align*}
		
	\noindent
	where $t(\rho) \vdash d$ is the cycle type of $\rho$, and likewise for $t(\sigma)$.
	The internal sum may be written 
	
		\begin{align*}
			\sum_{\rho,\sigma \in S(d)}  p_{t(\sigma)}(A_N) p_{t(\rho)}(B_N) \vec{w}^r(\rho,\sigma)
			&= \sum_{\alpha \vdash d} \sum_{\beta \vdash d} p_\alpha(A_N) p_\beta(B_N) \sum_{\rho \in C_\alpha} \sum_{\sigma \in C_\beta}
				\vec{w}^r(\rho,\sigma)\\
			&= \sum_{\alpha \vdash d} \sum_{\beta \vdash d} p_\alpha(A_N) p_\beta(B_N) \vec{W}^r(\alpha,\beta),
		\end{align*}

	\noindent		
	and this completes the proof of Theorem \ref{thm:LeadingDerivatives}.
	
\section{(Monotone) Hurwitz theory and the HCIZ free energy}
As shown by Hurwitz in the 19th century
\cite{Hurwitz}, the enumeration of unrestricted walks on the symmetric group with
given boundary conditions is equivalent to the enumeration of branched covers
of the Riemann sphere with given singular data.

To state this precisely, consider the generating function 
	
	\begin{equation*}
		\mathbf{W} = 1 + \sum_{d=1}^{\infty} \frac{z^d}{d!} \sum_{r=0}^{\infty} \frac{t^r}{r!}
		\sum_{\alpha,\beta \vdash d} p_\alpha(A)p_\beta(B) W^r(\alpha,\beta),
	\end{equation*}
	
\noindent
where

	\begin{equation*}
		A = \begin{bmatrix}
			a_1 & {} & {} \\
			{} & a_2 & {} \\
			{} & {} & \ddots
		\end{bmatrix}, \quad
		B = \begin{bmatrix}
			b_1 & {} & {} \\
			{} & b_2 & {} \\
			{} & {} & \ddots
		\end{bmatrix}
	\end{equation*}
	
\noindent
are a pair of formal infinite diagonal matrices and $W^r(\alpha,\beta)$ is the total 
number of $r$-step walks on 
$S(d)$ which begin in $C_\alpha$ and end in $C_\beta$.  Thus $\mathbf{W}$ is an
element of the formal power series algebra $\field{Q}[[z,t,a_1,a_2,\dots,b_1,b_2,\dots]]$.
Set

	\begin{equation*}
		\hur^r(\alpha,\beta) = \bigg{[} \frac{z^d}{d!} \frac{t^r}{r!} p_\alpha(A) p_\beta(B)\bigg{]}\Hur,
	\end{equation*}
	
\noindent
where $\Hur = \log \mathbf{W}$ and $[X]Y$ denotes the coefficient of term $X$
in a series $Y$.  By the exponential formula \cite{GJ:book}, the coefficient $\hur^r(\alpha,\beta)$
is the number of $r$-step walks beginning in $C_\alpha$ and 
ending in $C_\beta$ whose endpoints and steps together generate a 
transitive subgroup of $S(d)$.  

Hurwitz showed that
$\hur^r(\alpha,\beta)/d!$ may be interpreted as a weighted count of degree $d$ branched 
covers of the 
Riemann sphere by a compact, connected Riemann surface such that the covering
map has profile $\alpha$ over $\infty$, $\beta$ over $0$, and simple branching over
each of the $r$th roots of unity.  According to the Riemann-Hurwitz formula, such a cover
exists if and only if 

	\begin{equation*}
		g = \frac{r+2-\ell(\alpha)-\ell(\beta)}{2}
	\end{equation*}
	
\noindent
is a nonnegative integer, in which case it is the topological genus of the covering surface.
We will use the notation $\hur^r(\alpha,\beta)=\hur_g(\alpha,\beta)$, with
the understanding that $r$ and $g$ determine one another via the Riemann-Hurwitz
formula.

The numbers $\hur_g(\alpha,\beta)$ 
were first considered from a modern perspective by Okounkov \cite{Okounkov}, who
called them the \emph{double Hurwitz numbers}.
Proving a conjecture of Pandharipande \cite{Pandharipande} in Gromov-Witten 
theory, Okounkov showed that $\Hur$ is a solution of the 2D Toda lattice hierarchy of Ueno and
Takasaki.  It was subsequently shown by Kazarian and Lando \cite{KL} that, when 
combined with the ELSV formula \cite{ELSV}, Okounkov's result implies the celebrated
Kontsevich-Witten theorem relating intersection theory in moduli spaces of Riemann
surfaces to the KdV hierarchy.

\subsection{Monotone Hurwitz numbers}
Mimicking the above classical construction, consider the generating function

	\begin{equation*}
		\vec{\mathbf{W}} = 1 + \sum_{d=1}^{\infty} \frac{z^d}{d!} \sum_{r=0}^\infty t^r 
		\sum_{\alpha,\beta \vdash d} p_\alpha(A_N) p_\beta(B_N) \vec{W}^r(\alpha,\beta)
	\end{equation*}
	
\noindent
enumerating \emph{monotone} walks of all possible lengths and boundary conditions 
on all of the finite symmetric groups.  Note that, due to the monotonicity
constraint, the variable $t$ is now an ordinary rather than exponential marker for the 
walk length statistic.
Define the \emph{monotone double Hurwitz numbers} by

	\begin{equation*}
		\mon^r(\alpha,\beta) = \bigg{[} \frac{z^d}{d!} t^r p_\alpha(A) p_\beta(B) \bigg{]} \Mon,
	\end{equation*}
	
\noindent
where $\Mon = \log \vec{\mathbf{W}}$.  Then $\mon^r(\alpha,\beta)$ is the number of $r$-step
monotone walks beginning in $C_\alpha$ and ending in $C_\beta$ whose endpoints and 
steps together generate a transitive subgroup of $S(d)$.  Alternatively, 
$\mon^r(\alpha,\beta)/d!$ counts a combinatorially
restricted subset of the branched covers counted by $\hur^r(\alpha,\beta)/d!$.  
In \cite{GGN1,GGN2}, an extensive combinatorial theory of monotone Hurwitz
numbers was developed, and monotone analogues of most combinatorial
results in classical Hurwitz theory were obtained.  Here we will use those
results in tandem with Theorem \ref{thm:LeadingDerivatives}.

\subsection{The HCIZ free energy}
The monotone double Hurwitz numbers describe the Maclaurin coefficients of the
HCIZ free energy.

As in the Introduction,
let $\Omega_N$ be a simply connected open set in $\field{C}$ which contains
the origin and avoids the zeros of $I_N(z)$, and let $\log I_N(z)$ denote the principal
branch of the logarithm of $I_N(z)$ on $\Omega_N$.
From Theorem \ref{thm:LeadingDerivatives} and the exponential formula, it
follows that the Maclaurin series of the logarithm is

	\begin{equation*}
		\log I_N(z) = \sum_{d=1}^N \frac{z^d}{d!} \sum_{r=0}^\infty \bigg{(} -\frac{1}{N} \bigg{)}^r
		\sum_{\alpha,\beta \vdash d} p_\alpha(A_N) p_\beta(B_N) \mon^r(\alpha,\beta) + O(z^{N+1}),
	\end{equation*}
	
\noindent
where the $O$-term is uniform on compact subsets of $\Omega_N$.

The Maclaurin series of $\log I_N(z)$ can be simplified using the Riemann-Hurwitz formula.
We have:

	\begin{align*}
		\log I_N(z) &=  \sum_{d=1}^N \frac{z^d}{d!} \sum_{\alpha,\beta \vdash d} p_\alpha(A_N) p_\beta(B_N) \sum_{r=0}^\infty \bigg{(} -\frac{1}{N} \bigg{)}^r
		\mon^r(\alpha,\beta) + O(z^{N+1}) \\
		&=\sum_{d=1}^N \frac{z^d}{d!} \sum_{\alpha,\beta \vdash d} p_\alpha(A_N) p_\beta(B_N) \sum_{g=0}^\infty \bigg{(} -\frac{1}{N} \bigg{)}^{2g-2+\ell(\alpha)+\ell(\beta)}
		\mon_g(\alpha,\beta) + O(z^{N+1}) \\
		&= N^2\sum_{d=1}^N \frac{z^d}{d!} 
		\sum_{g=0}^{\infty} \frac{1}{N^{2g}}\sum_{\alpha,\beta \vdash d}(-1)^{\ell(\alpha)+\ell(\beta)} \frac{p_\alpha(A_N)}{N^{\ell(\alpha)}} 
		\frac{p_\beta(B_N)}{N^{\ell(\beta)}} \mon_g(\alpha,\beta) + O(z^{N+1}),
	\end{align*}
	
\noindent
with the $O$-term uniform on compact subsets of $\Omega_N$.  Thus, 
parameterizing the monotone double Hurwitz numbers by $g$ instead of $r$,
we arrive at a description of $F_N(z)=N^{-2}\log I_N(z)$ which is
well-poised for an $N \rightarrow \infty$ asymptotic analysis.

\subsection{Asymptotic expansion of Maclaurin coefficients}
We will now consider the asymptotic behaviour of the Maclaurin coefficients of $F_N(z)$
as $N \rightarrow \infty$, when $A_N,B_N$ vary regularly with $N$.

Suppose there exists a nonnegative number $M$, a nonnegative integer $h$, 
and two complex sequences $\phi_k,\psi_k$ such that:

	\begin{enumerate}
	
		\smallskip
		\item
		$\|A_N\|,\|B_N\| \leq M$ for all $N \geq 1$;
		
		\smallskip
		\item
		For each $k \geq 1$,
		
			\begin{equation*}
				\frac{1}{N}\Tr(A_N^k) = \phi_k + o\bigg{(} \frac{1}{N^{2h}} \bigg{)}, \quad
				\frac{1}{N}\Tr(B_N^k) = \psi_k + o\bigg{(} \frac{1}{N^{2h}} \bigg{)}
			\end{equation*}
			
		\noindent
		as $N \rightarrow \infty$.
	\end{enumerate}
	
\noindent
We will summarize these conditions by saying that $A_N,B_N$ are $(M,h)$-regular 
with limit moments $\phi_k,\psi_k$.

\begin{theorem}
	\label{thm:MaclaurinAsymptotics}
	Suppose that $A_N,B_N$ are $(M,h)$-regular with limit moments $\phi_k,\psi_k$.
	Then, for each $d \geq 1$, we have
	
		\begin{equation*}
			F_N^{(d)}(0) = \sum_{g=0}^h \frac{C_{g,d}}{N^{2g}} + o\bigg{(} \frac{1}{N^{2h}}\bigg{)}
		\end{equation*}
		
	\noindent
	as $N \rightarrow \infty$, where
	
		\begin{equation*}
			C_{g,d} = \sum_{\alpha,\beta \vdash d} (-1)^{\ell(\alpha)+\ell(\beta)} \phi_\alpha \psi_\beta \mon_g(\alpha,\beta),
		\end{equation*}	
		
	\noindent
	and 
	
		\begin{equation*}
			\phi_\alpha = \prod_{i=1}^{\ell(\alpha)} \phi_{\alpha_i}, \quad \psi_\beta = \prod_{j=1}^{\ell(\beta)} \phi_{\beta_j}.
		\end{equation*}
\end{theorem}

\begin{proof}
	For $1 \leq d \leq N$, we have 
		
		\begin{equation*}
			F_N^{(d)}(0) = \sum_{g=0}^{\infty} \frac{C_{g,d,N}}{N^{2g}},
		\end{equation*}
		
	\noindent
	where 
	
		\begin{equation*}
			C_{g,d,N} = \sum_{\alpha,\beta \vdash d} (-1)^{\ell(\alpha)+\ell(\beta)}
			\frac{p_\alpha(A_N)}{N^{\ell(\alpha)}} \frac{p_\beta(B_N)}{N^{\ell(\beta)}} \mon_g(\alpha,\beta).
		\end{equation*}
		
	Since $\|A_N\|,\|B_N\| \leq M$ for all $N \geq 1$, we have
	
		\begin{equation*}
			|C_{g,d,N}| \leq M^{2d} \sum_{\alpha,\beta \vdash d} \mon_g(\alpha,\beta)
		\end{equation*}
		
	\noindent
	for all $N \geq 1$.  Since $\mon_g(\alpha,\beta)$ counts certain solutions of the equation
	
		\begin{equation*}
			\sigma = \rho \tau_1 \dots \tau_r 
		\end{equation*}
		
	\noindent
	in $S(d)$, with $r=2g-2+\ell(\alpha)+\ell(\beta)$, we have
	
		\begin{equation*}
			\mon_g(\alpha,\beta) \leq (d!)^{2g+2d}
		\end{equation*}
		
	\noindent
	for all $\alpha,\beta \vdash d$.  Consequently, 
	
		\begin{equation*}
			|C_{g,d,N}| \leq (d!p(d)M)^{2d} (d!)^{2g},
		\end{equation*}
		
	\noindent
	where $p(d)$ is the number of Young diagrams with $d$ cells.
	Moreover, since $A_N,B_N$ are $h$-regular with limit moments
	$\phi_k,\psi_k$, we have
	
		\begin{equation*}
			C_{g,d,N} = C_{g,d} + o\bigg{(} \frac{1}{N^{2h}}\bigg{)}
		\end{equation*}
		
	\noindent 
	as $N \rightarrow \infty$.
	We thus have
	
		\begin{align*}
			F_N^{(d)}(0) &= \sum_{g=0}^h \frac{C_{g,d,N}}{N^{2g}} + \sum_{g=h+1}^{\infty} \frac{C_{g,d,N}}{N^{2g}} \\
			&= \sum_{g=0}^h \frac{C_{g,d}}{N^{2g}} + o\bigg{(} \frac{1}{N^{2h}} \bigg{)} + O\bigg{(} \frac{1}{N^{2h+2}}\bigg{)}
		\end{align*}
		
	\noindent
	as $N \rightarrow \infty$, which proves the claim.
\end{proof}

\begin{remark}
The convergence of $F_N^{(d)}(0)$ under the above hypotheses was first stated by Itzykson and Zuber \cite{IZ},
and proved by Collins \cite{Collins}.  Collins obtained the limit of $F_N^{(d)}(0)$ as a double sum over
$S(d)$, whereas we present the same limit as a double sum over partitions of $d$. 
\end{remark}

\begin{remark}
To the best of our knowledge, Theorem \ref{thm:MaclaurinAsymptotics} is the
first result which addresses the sub-leading asymptotics of $F_N^{(d)}(0)$, clearly showing the emergence of
a ``topological expansion'' in this context.
\end{remark}

\section{Absolute summability of genus-specific generating functions}
Let us introduce a sequence of formal power series $C_g(z) \in \field{C}[[z]]$ defined
by 

	\begin{equation}
		\label{eqn:GenusSpecificAnalytic}
		C_g(z) = \sum_{d=1}^{\infty} \frac{z^d}{d!}C_{g,d} , \quad g \geq 0,
	\end{equation}
	
\noindent
where $C_{g,d}$ are the expansion coefficients from Theorem \ref{thm:MaclaurinAsymptotics}.
The goal of this section is to prove that these formal power series are absolutely summable, 
and to obtain bounds on their radii of convergence. By the uniform boundedness of $\|A_N\|$
and $\|B_N\|$, this reduces to establishing the absolute summability of the series 

	\begin{equation}
		\label{eqn:GenusSpecificAlgebraic}
		\Mon_g(z) = \sum_{d=1}^{\infty} \frac{z^d}{d!}\sum_{\alpha,\beta \vdash d} \mon_g(\alpha,\beta) , \quad g \geq 0.
	\end{equation}

\subsection{Monotone simple Hurwitz numbers}
Consider the genus-specific generating functions 

	\begin{equation*}
		\vec{\mathbf{S}}_g = \sum_{d=1}^{\infty} \frac{z^d}{d!} \mon_{g,d}
	\end{equation*}
	
\noindent
for the monotone \emph{simple} Hurwitz numbers

	\begin{equation*}
		\mon_{g,d} = \mon_g(1^d,1^d).
	\end{equation*}
	
\noindent
The monotone simple Hurwitz number $\mon_{g,d}$ counts
monotone loops of length $r=2g-2+2d$, based at any given
point of $S(d)$, whose steps generate a transitive subgroup
of $S(d)$.

According to \cite[Theorem 1.4]{GGN2}, for any
$g \geq 2$ we have 

	\begin{equation*}
		\vec{\mathbf{S}}_g = \frac{\zeta(1-2g)}{2-2g} + \frac{1}{(1-6s)^{2g-2}}
		\sum_{r=0}^{3g-3} \sum_{\mu \vdash r} \frac{c_{g,\mu}(6s)^{\ell(\mu)}}{(1-6s)^{\ell(\mu)}},
	\end{equation*}
	
\noindent
where $\zeta$ is the Riemann zeta function, the $c_{g,\mu}$'s are rational numbers, and
$s$ is the unique solution of the functional equation

	\begin{equation*}
		s=z(1-2s)^{-2}
	\end{equation*}
	
\noindent
in the formal power series algebra $\field{Q}[[z]]$. 
This equation may be solved by Lagrange inversion,
yielding the solution
	
	\begin{equation*}
		s = \sum_{n=1}^{\infty} \frac{2^{n-1}}{n} {3n-2 \choose n-1}z^n.
	\end{equation*}
	
\noindent 
Since 

	\begin{equation*}
		s' = \sum_{n=0}^{\infty} {3n+1 \choose n} (2z)^n = {}_2F_1\bigg{(} \frac{2}{3}, \frac{4}{3}, \frac{3}{2}; \frac{27}{2}z\bigg{)},
	\end{equation*}
	
\noindent
where ${}_2F_1(a,b,c;z)$ is the Gauss hypergeometric function, $\vec{\mathbf{S}}_g$ extends to a holomorphic function
of $z$ on the domain $\field{C}\backslash [z_c,\infty)$, where 

	\begin{equation*}
		z_c = \frac{2}{27}.
	\end{equation*}
	
In the case $g=0$, \cite[Theorem 1.1]{GGN1} yields the exact formula

	\begin{equation*}
		\frac{\mon_{0,d}}{d!} = \frac{2^{d-1}}{d^2(d-1)}{3d-3 \choose d-1},
	\end{equation*}
	
\noindent
so that, using Stirling's formula, we obtain $\frac{2}{27}$ as the 
radius of convergence of $\vec{\mathbf{S}}_0$.  Moreover, 
since $\vec{\mathbf{S}}_0$ has positive coefficients, 
Pringsheim's theorem (see \cite[\S7.21]{Titchmarsh}) 
guarantees that $z_c=\frac{2}{27}$ is a singularity of $\vec{\mathbf{S}}_0$.

 We have thus shown that:

	\begin{theorem}
		\label{thm:simple}
		The series $\vec{\mathbf{S}}_g$, $g \geq 0$, have a common dominant 
		singularity at the critical point $z_c=\frac{2}{27}$.
	\end{theorem}
	
\begin{remark}
	The $g=0$ case of Theorem \ref{thm:simple} was obtained by
	Zinn-Justin in \cite{ZJ} using the Toda lattice equations.  For
	a combinatorial solution of the Toda equations encompassing
	those of Okounkov and Zinn-Justin, see \cite{GGN3}.
\end{remark}

\begin{remark}
	P. Di Francesco has pointed out to us that the critical value 
	$z_c=\frac{2}{27}$ also appears in the enumeration of finite
	groups, see \cite{Pyber}.
\end{remark}

\subsection{Monotone double Hurwitz numbers}
Consider the full genus-specific generating functions
$\Mon_g(z)$ defined in equation \ref{eqn:GenusSpecificAlgebraic}.
Obviously,

	\begin{equation*}
		\sum_{\alpha,\beta \vdash d} \mon_g(\alpha,\beta) \geq \mon_{g,d},
	\end{equation*}
	
\noindent
so the radius of convergence of $\Mon_g$ is at most the radius of convergence
of $\vec{\mathbf{S}}_g$.  

We now explain how a lower bound on the radius of convergence of $\Mon_g$
follows from the study of a refinement of the monotone double Hurwitz number 
$\mon_g(\alpha,\beta)$.  Given a positive integer $c$, let $\mon_g(\alpha,\beta;c)$
denote the number of walks counted by $\mon_g(\alpha,\beta)$ whose steps 
have $c$ distinct labels.  The following inequality is obtained in \cite{GN}:

	\begin{equation*}
		\sum_{c=2}^d 3^c  \mon_g(1^d,1^d;c) \leq \sum_{\alpha,\beta \vdash d} \mon_g(\alpha,\beta) \leq \sum_{c=2}^d 4^c \mon_g(1^d,1^d;c).
	\end{equation*}
	
\noindent
The proof is combinatorial, and makes use of an action of the symmetric
group $S(r)$ on the set of walks counted by the classical 
double Hurwitz number $\hur_g(\alpha,\beta)$.  From
this inequality and the definition of $\mon_g(\alpha,\beta;c)$, we obtain 

	\begin{equation*}
		\sum_{\alpha,\beta \vdash d} \mon_g(\alpha,\beta) \leq 4^{d-1} \sum_{c=2}^d  \mon_g(1^d,1^d;c)
		= 4^{d-1} \mon_{g,d},
	\end{equation*}
	
\noindent
which implies that the radius of convergence of $\Mon_g$ is at least 
one quarter the radius of convergence of $\vec{\mathbf{S}}_g$.
Combining this with Theorem \ref{thm:simple}, we have:

	\begin{theorem}	
		\label{thm:RadiusOfConvergence}
		For each $g \geq 0$, the series $\Mon_g$ is absolutely
		summable, and has radius of convergence at least 
		$\frac{1}{54}$ and at most $\frac{2}{27}$.
	\end{theorem}	
	
\begin{remark}
	The results of \cite{CGM} also imply the absolute summability
	of $\Mon_0$, but without effective bounds on the radius of
	convergence.
\end{remark}

\section{Convergence of the HCIZ free energy}

\subsection{Convergence under a non-vanishing hypothesis}
Let $A_N,B_N$ be $(M,0)$-regular with limit moments $\phi_k,\psi_k$,
respectively.  Then, by Theorem \ref{thm:RadiusOfConvergence}, the
series $C_0(z)$ is absolutely summable, with radius of convergence 
at least $\frac{1}{54M^2}$.
	
\begin{theorem}
	\label{thm:Convergence}
	Suppose that there exists of a positive number $R$ such that $I_N(z)$ is non-vanishing
	on the open disc $D(0,R)$ for all but finitely many $N$.  Let 

	\begin{equation*}
		r = \min \bigg{\{} R, \frac{1}{54M^2} \bigg{\}}.
	\end{equation*}
	
	\noindent
	Then $F_N(z) \rightarrow C_0(z)$ uniformly on compact subsets of $D(0,r)$.
\end{theorem}

\begin{proof}
	First, we show that the claim holds pointwise.  
	The proof is a two circles argument.	
	Let $z_0 \in D(0,r)$ and $\varepsilon>0$ be arbitrary.
	Choose $r_1,r_2$ such that
	
		\begin{equation*}
			|z_0| < r_1 < r_2 < r.
		\end{equation*}
		
	\noindent
	For $N$ sufficiently large, we have
	
		\begin{equation*}
			|F_N(z_0) - C_0(z_0)| \leq \sum_{d=1}^{\infty} |F_N^{(d)}(0) - C_{0,d}| \frac{|z_0|^d}{d!}.
		\end{equation*}
		
	\noindent
	By Cauchy's inequality, 
	
		\begin{equation*}
			\frac{1}{d!} |F_N^{(d)}(0) - C_{0,d}| \leq \frac{\|F_N-C_0\|_{r_1}}{r_1^d},
		\end{equation*}
		
	\noindent
	where $\| \cdot \|_{r_1}$ denotes sup-norm over the circle of radius $r_1$.  
	We thus have
	
		\begin{equation*}
			|F_N(z_0) - C_0(z_0)| \leq \sum_{d=1}^E |F_N^{(d)}(0) - C_{0,d}| \frac{|z_0|^d}{d!}
			+ \frac{\|F_N\|_{r_1} + \|C_0\|_{r_1}}{1-\frac{|z_0|}{r_1}}\bigg{(} \frac{|z_0|}{r_1} \bigg{)}^{E+1}
		\end{equation*}
		
	\noindent
	for any positive integer $E$.
	
	We will now bound the error term independently of $N$.  Since
	
		\begin{align*}
			|I_N(z)| &\leq \int\limits_{U(N)} e^{|z| N|\Tr(A_NUB_NU^{-1})|} \mathrm{d}U \\
			& \leq \int\limits_{U(N)} e^{|z| M^2 N \sum_{i,j=1}^N |u_{ij}|^2} \mathrm{d}U \\
			& \leq e^{M^2N^2|z|}
		\end{align*}
		
	\noindent
	for all $z \in \field{C}$, where the last line follows from the fact that $(|u_{ij}|^2)$ is a
	doubly stochastic matrix, the inequality 
	
		\begin{equation*}
			\Re F_N(z) \leq M^2|z| 
		\end{equation*}
		
	\noindent
	holds for all $z \in \Omega_N$, the domain of holomorphy of $F_N(z)$.  Combining
	this with the Borel-Carath\'eodory inequality, which bounds the sup-norm of
	an analytic function on a circle in terms of the maximum of its real part on a
	circle of larger radius (see e.g. \cite[\S 5.5]{Titchmarsh}, we have
	
		\begin{equation*}
			\|F_N\|_{r_1} \leq \frac{2r_1}{r_2-r_1} \sup_{|z|=r_2} \Re F_N(z) \leq  \frac{2M^2r_1r_2}{r_2-r_1}.
		\end{equation*}
		
	Returning to our estimate on $|F_N(z_0)-C_0(z_0)|$, we now have the inequality
	
		\begin{equation*}
			|F_N(z_0) - C_0(z_0)| \leq \sum_{d=1}^E |F_N^{(d)}(0) - C_{0,d}| \frac{|z_0|^d}{d!}
			+ \frac{\frac{2M^2r_1r_2}{r_2-r_1} + \|C_0\|_{r_1}}{1-\frac{|z_0|}{r_1}}\bigg{(} \frac{|z_0|}{r_1} \bigg{)}^{E+1}
		\end{equation*}
		
	\noindent 
	for all $N$ sufficiently large and all $E \geq 1$.  Choosing $E_0$ large enough so that 
	
		\begin{equation*}
			\frac{\frac{2M^2r_1r_2}{r_2-r_1} + \|C_0\|_{r_1}}{1-\frac{|z_0|}{r_1}}\bigg{(} \frac{|z_0|}{r_1} \bigg{)}^{E+1} < \frac{\varepsilon}{2},
		\end{equation*}
		
	\noindent
	and subsequently choosing $N_0$ large enough so that
	
		\begin{equation*}
			 \sum_{d=1}^{E_0} |F_{N_0}^{(d)}(0) - C_{0,d}| \frac{|z_0|^d}{d!} < \frac{\varepsilon}{2},
		\end{equation*}
		
	\noindent
	we obtain that 
	
		\begin{equation*}
			|F_N(z_0) - C_0(z_0)| < \varepsilon
		\end{equation*}
		
	\noindent
	for all $N \geq N_0$.  This completes the proof.
	
	We now explain how the mode of convergence
	can be boosted \ from pointwise to uniform on compact subsets of $D(0,r)$.  It is
	easy to check that

	\begin{equation*}
		\sup_N \|F_N\|_K < \infty
	\end{equation*}
	
\noindent
for each compact set $K \subset D(0,r)$; for example, one could prove
this holds for closed discs by using the Borel-Carath\'edory inequality
again.  Thus $\{F_N\}$ is a locally uniformly bounded family
on $D(0,r)$, and Vitali's theorem \cite[\S 5.21]{Titchmarsh} implies
pointwise and uniform-on-compact convergence are equivalent.
\end{proof}

\subsection{Remarks on the non-vanishing hypothesis}
Theorem \ref{thm:Convergence} reduces the complex convergence 
of $F_N(z)$ near $z=0$ to checking that the zeros of $I_N(z)$
do not encroach on $z=0$ as $N \rightarrow \infty$.  
Therefore it is of significant interest
to determine sufficient hypotheses on the matrix sequences $A_N$
and $B_N$ which ensure that this condition holds.

Determining the exact locations of the zeros of $I_N(z)$ seems
to be a difficult problem in general.  Let us give one example
where it can be solved explicitly.  Suppose that $A_N$
and $B_N$ are real diagonal with distinct eigenvalues 

	\begin{equation*}
		a_1 < \dots < a_N, \quad b_1< \dots < b_N,
	\end{equation*}	
	
\noindent
and suppose further that the eigenvalues of $A_N$ form
an arithmetic progression with constant difference $\hbar >0$.
In this case, all determinants in the HCIZ formula (see \cite[Equation 3.28]{IZ})
are Vandermonde determinants, and one computes 
that the zeros of $I_N(z)$ are the pure imaginary points

	\begin{equation*}
		z = \frac{2k\pi}{N\hbar(b_j-b_i)}\mathbf{i}, \quad 1 \leq i<j \leq N, \quad k \in \field{Z}.
	\end{equation*}
	
\noindent
Thus one has a sort of Lee-Yang theorem in this special case.
If the eigenvalues of $A_N$ and $B_N$ are confined to a fixed 
interval $[-M,M]$ for all $N$, then $\hbar$ is of order $1/N$ and 
the above formula shows that the zeros
of $I_N(z)$ remain at bounded distance from the origin as $N$
increases.

Concerning the verification of the non-vanishing hypothesis
more generally, we do not know much at present.  It might 
be that the Toda equations could be of help here, but we have
not pursued this seriously.

\subsection{A conjecture}
In this article, we have produced explicit holomorphic candidates $C_g(z)$
for the asymptotics of $F_N(z)$ near $z=0$.  Indeed, if $F_N(z)$ admits a
uniform-on-compact asymptotic expansion on the asymptotic scale $N^{-2}$
in some neighbourhood $\Omega$ of $z=0$, then the coefficients in this 
expansion \emph{must} be the functions $C_g(z)$.  We therefore believe
that the following conjecture is reasonable.

	\begin{conjecture}
		If $A_N,B_N$ are $(M,h)$-regular with limit moments $\phi_k,\psi_k$,
		then there exists a neighbourhood $\Omega$ of $z=0$ such that
		
			\begin{equation*}
				F_N(z) = \sum_{g=0}^h \frac{C_g(z)}{N^{2g}} + o\left( \frac{1}{N^{2h}} \right)
			\end{equation*}
			
		\noindent
		uniformly on compact subsets of $\Omega$.
	\end{conjecture}

\bibliographystyle{amsplain}

\begin{thebibliography}{10}




























\bibitem{Collins}
B. Collins,
\textit{Moments and cumulants of polynomial random variables on unitary groups, the Itzykson-Zuber integral,
and free probability},
Int. Math. Res. Not. IMRN. \textbf{17} (2003), 954-982.

\bibitem{CGM}
B. Collins, A. Guionnet, E. Maurel-Segala, \textit{Asymptotics of unitary and orthogonal matrix integrals},
Adv. Math. \textbf{222} (2009), 172-215.


\bibitem{CS}
B. Collins, P. \'Sniady,
\textit{Integration with respect to the Haar measure on unitary, orthogonal and symplectic group},
Comm. Math. Phys, \textbf{264} (2006), 773-795.















\bibitem{ELSV}
T. Ekedahl, S. Lando, M. Shapiro, A. Vainshtein,
\textit{Hurwitz numbers and intersections on moduli spaces of curves}, 
Invent. Math. \textbf{146} (2001), 297-327.






\bibitem{EY}
L. Erd\H{o}s, H.-T. Yau,
\textit{Universality of local spectral statistics of random matrices},
Bulletin (new series) of the American Mathematical Society \textbf{49}(3) (2012),
377-414.








\bibitem{GGN1}
I. P. Goulden, M. Guay-Paquet, J. Novak,
\textit{Monotone Hurwitz numbers in genus zero},
Canad. J. Math. \textbf{65} (2013), 1020-1042.

\bibitem{GGN2}
I. P. Goulden, M. Guay-Paquet, J. Novak,
\textit{Polynomiality of monotone Hurwitz numbers in higher genera},
Adv. Math. \textbf{238} (2013), 1-23.

\bibitem{GGN3}
I. P. Goulden, M. Guay-Paquet, J. Novak,
\textit{Toda equations and piecewise polynomiality for mixed double Hurwitz numbers},
submitted.

\bibitem{GJ:book}
I. P. Goulden, D. M. Jackson,
\textit{Combinatorial Enumeration}, John Wiley and Sons, New York, 1983 
(reprinted by Dover, 2004).










\bibitem{GN}
M. Guay-Paquet, J. Novak,
\textit{A self-interacting random walk on the symmetric group},
in preparation.

\bibitem{Guionnet:CMP}
A. Guionnet,
\textit{First order asymptotics of matrix integrals; a rigorous approach towards
the understanding of matrix models},
Communications in Mathematical Physics \textbf{244} (2004), 527-569.

\bibitem{Guionnet:PS}
A. Guionnet,
\textit{Large deviations and stochastic calculus for large random matrices},
Probability Surveys \textbf{1} (2004), 72-172.









\bibitem{HC}
Harish-Chandra,
\textit{Differential operators on a semisimple Lie algebra},
American Journal of Mathematics \textbf{79} (1957), 87-120.



\bibitem{Hurwitz}
A. Hurwitz,
\textit{Ueber Riemann'sche Fl\"achen mit gegebenen Verzweigegungspunkten},
Mathematische Annalen \textbf{39} (1891), 1-60.


\bibitem{IZ}
C. Itzykson, J.-B. Zuber,
\textit{The planar approximation. II},
Journal of Mathematical Physics \textbf{21}(3) (1980), 411-421.









\bibitem{KL}
M. E. Kazarian, S. K. Lando,
\textit{An algebro-geometric proof of Witten's conjecture},
Journal of the AMS \textbf{20}(4) (2007), 1079-1089.












\bibitem{MN}
S. Matsumoto, J. Novak,
\textit{Jucys-Murphy elements and unitary matrix integrals},
Int. Math. Res. Not. IMRN \textbf{2} (2013), 362-297.










\bibitem{Novak:BCP}
J. Novak,
\textit{Jucys-Murphy elements and the unitary Weingarten function},
Banach Center Publications \textbf{89} (2010), 231-235.




\bibitem{Okounkov}
A. Okounkov,
\textit{Toda equations for Hurwitz numbers},
Mathematical Research Letters \textbf{7} (2000), 447-453.










\bibitem{Pandharipande}
R. Pandharipande,
\textit{The Toda equations and the Gromov-Witten theory of the Riemann sphere},
Letters in Mathematical Physics \textbf{53} (2000), 59-74.




\bibitem{Pyber}
L. Pyber,
\textit{Enumerating finite groups of given order},
Annals of Mathematics \textbf{137}(1) (1993), 203-220.
















\bibitem{Titchmarsh}
E. C. Titchmarsh,
\textit{The Theory of Functions, Second Edition},
Oxford University Press, 1939.










\bibitem{ZJ}
P. Zinn-Justin,
\textit{HCIZ integral and 2D Toda lattice hierarchy},
Nuclear Physics B \textbf{634} [FS] (2002), 417-432.


\bibitem{ZZ}
P. Zinn-Justin, J.-B. Zuber,
\textit{On some integrals over the $U(N)$ unitary group and their large $N$ limit},
Journal of Physics A: Mathematical and General \textbf{36} (2003), 3173-3193.



\end{thebibliography}

\end{document}